\documentclass[12pt]{amsart}

\usepackage{amsmath,amssymb}

\DeclareMathOperator{\Qf}{Qf}

\DeclareMathOperator{\win}{\mathit w_{in}}

\newtheorem{thm}{Theorem}[section]

\newtheorem{prop}[thm]{Proposition}

\newtheorem{lem}[thm]{Lemma}

\newtheorem{cor}[thm]{Corollary}

\newcommand\qf{\operatorname{Qf}}

\newcommand\mtwo{\mathcal M_2}

\newcommand\mone{\mathcal M_1}

\begin{document}


    \title[Complete integral closure]
{Complete integral closure \\ and strongly divisorial prime ideals}



\author{Valentina Barucci}

\address{Dipartimento di Matematica, Universit\`{a} di Roma "La 
Sapienza",
Piazzale A.  Moro, 5, 00185 Roma, Italy}

\email{barucci@mat.uniroma1.it}

\author{Stefania Gabelli}

\address{Dipartimento di Matematica, Universit\`{a} degli Studi Roma 
Tre,
Largo S.  L.  Murialdo,
1, 00146 Roma, Italy}

\email{gabelli@mat.uniroma3.it}

\author{Moshe Roitman}

\address{Department of Mathematics, University of
Haifa, Mount Carmel, Haifa 31999, Israel}

\email{mroitman@math.haifa.ac.il}




\begin{abstract}

It is well known that a domain without proper
strongly divisorial ideals is completely integrally closed.
In this paper we show that a domain without {\em prime} strongly
divisorial ideals is not necessarily completely integrally closed,
although this property holds
under some additional assumptions.
    \end{abstract}

    \maketitle




\section*{Introduction}

The strongly divisorial ideals of a domain $R$ are the nonzero
conductors $(R:T)$ for some overring $T$ of $R$.  They can be used to
describe the complete integral closure $R^*$ of $R$; in fact
$R^* =$
$\bigcup(R:I)$, where $I$ runs throughout the set of the strongly divisorial
ideals of $R$.  In particular, $R$ is completely integrally closed,
(that is $R = R^*$), if and only if $R$ has no proper strongly
divisorial ideals \cite{B}.  This paper deals with the following
question \cite[Problem 11, page 461]{CG}:

\smallskip
{\bf Question}. Is a domain $R$ completely integrally closed if $R$
does not have strongly divisorial {\em prime} ideals?

\smallskip
We show that in general the answer is negative, and give conditions
on $R$ for a positive answer.

The positive results are collected in Section 1.  We first prove
some criteria useful
to establish when a domain having a proper strongly divisorial ideal
is forced to have a strongly divisorial prime ideal.
We prove that an ideal which is maximal in the set of proper strongly
divisorial
ideals of a domain is prime (Proposition \ref{maxstrdiv}) and we
give conditions sufficient for a strongly divisorial ideal to
have a strongly divisorial minimal prime (Lemmas \ref{preirred} and
\ref{irred}).

Applying these results, we show that the answer to the question above is
positive for some classes of domains.
We prove in Proposition \ref{forte} that if a domain has a proper strongly
divisorial ideal containing a power of its radical and if this radical is an
irredundant intersection of primes, then the domain contains a strongly
divisorial prime ideal. We note in Lemma \ref{rad} that in a seminormal domain each
strongly divisorial ideal is radical, and so the first assumption of
Proposition \ref{forte} is always satisfied in the seminormal case.
As another consequence of Proposition \ref{forte}, we give a positive answer
for
any pullback of a Noetherian domain (Theorem \ref{noeth}).

We also answer the above question affirmatively for the class of domains
with the
property that each maximal $t$-ideal is divisorial (Proposition
\ref{maxdiv}).  This class includes Noetherian,
Mori, TV- and pseudo-valuation domains.  Indeed a domain of this kind
without strongly divisorial prime ideals is a Krull domain, as it is
pointed out in Theorem \ref{Krull}, where there are given several
other equivalent conditions in terms of ideals.  A domain without
strongly divisorial primes satisfying the stronger property that each
maximal ideal is invertible is a Dedekind domain (Corollary
\ref{Dedekind}).

There are also other important classes of domains for which the question
has a positive answer; for example valuation domains (more generally,
TP-domains or Pr\"ufer RTP-domains) and one-dimensional seminormal
domains with finitely many maximal ideals (Proposition
\ref{seminormal}).

On the other hand we do not know the answer in case of a one-dimensional
seminormal domain.

In Section 2 we give examples of domains that are not completely
integrally closed but do not have strongly divisorial prime ideals. The
first two examples are quasilocal and one-dimensional domains whose
maximal ideal is not divisorial. The conductor of the complete integral
closure is zero in the first example and is nonzero in the second one.
These two domains are Ar\-chi\-medean and not seminormal by Proposition
\ref{seminormal}. The third example is an infinite dimensional domain
with nonzero conductor in its complete integral closure, and can be
chosen to satisfy or not either one of the conditions of being
Archimedean or seminormal. In the seminormal case, the conductor is a
strongly divisorial radical ideal (Lemma \ref{rad}) such that all its
minimal primes are strong but not divisorial.

\smallskip Throughout this paper, $R$ will denote an integral domain
with quotient field $\qf(R) = K$. For a pair of fractional ideals $I$
and $J$ of a domain $R$ we let $(J:I)$ denote the set $\{x\in K|\,
xI\subseteq J\}$ and $(J:_{R}I)$ denote the set $\{x\in R|\, xI\subseteq
J\}$. As usual, we set $I_v=(R:(R:I))$ and $I_t=\bigcup J_v$ with the
union taken over all finitely generated fractional ideals $J$ contained
in $I$. An ideal $I$ is a {\it $v$-ideal}, or {\it divisorial}, if
$I=I_v$, and it is a {\it $t$-ideal} if $I=I_t$. The $v$- and the
$t$-operations are particular $*$-operations. If $(R:I)=(I:I)$, that is,
$I = I(R:I)$, we say that $I$ is a {\it strong} ideal or a {\it trace}
ideal. A {\it strongly divisorial} ideal is an ideal that is strong and
divisorial and a {\it maximal strongly divisorial} ideal is an ideal
that is maximal in the set of the proper strongly divisorial ideals. A
{\it maximal $t$-ideal} is an ideal that is maximal in the set of the
proper integral $t$-ideals. A general reference for systems of ideals
and $*$-operations is \cite{HK}.

A nonzero element $x \in K$ is said to be {\it almost integral} over $R$ if
there is a nonzero element $d \in R$ such that $dx^i \in R$, for $i \ge 0$.
This is equivalent to saying that all the powers of $x$ belong to a finite
$R$-module or that they generate a fractional ideal of $R$.  The {\it
complete integral closure} of $R$ in $K$, which we denote by $R^{*}$, is
the set of the elements of $K$ that are almost integral over $R$.  One says
that $R$ is {\it completely integrally closed} if $R = R^{*}$.

If an element is integral over $R$, then it is almost integral, and if $R$
is Noetherian the converse holds.  Thus, if $R'$ denotes the integral
closure of $R$ in $K$, one has $R \subseteq R' \subseteq R^{*}$, and if $R$
is Noetherian then $R' = R^{*}$.  It is well known that $R^{*}$ is always
integrally closed, but many examples have been given to show that it might
not be completely integrally closed.  However, $R^{*}$ is completely
integrally closed if the conductor $(R:R^{*})$ is a nonzero ideal, that
is,
if $R^{*}$ is contained in a finite $R$-submodule of $K$ (see for example
\cite{B}).




\section{Positive results}

In this Section we give some criteria which are useful
for establishing when a domain having a proper strongly divisorial ideal
has a strongly divisorial prime ideal and we apply our criteria to give
a  positive answer to the Question in the Introduction for some
classes of domains.

We start by proving that a maximal strongly divisorial
ideal is prime.

\begin{lem} \label{genprime}
Let $R$ be any ring and let $\mathbf S$ be a set of proper ideals of $R$.
Assume that $(I:_Rx)\in \mathbf S$ for any $I\in \mathbf S$ and nonzero
$x\in R\setminus I$.  Then an ideal of $R$ that is maximal in $\mathbf S$
is prime.  
\end{lem}

\begin{proof} Let $I$ be an ideal that is maximal in $\mathbf S$ and let
$x,y\in R$ such that $xy\in I$ and $x\notin I$.  We have $I\subseteq
(I:_Rx)\in \mathbf S$.  Hence $(I:_R x)=I$ and $y\in (I:_Rx)=I$.  It
follows that $I$ is a prime ideal.  \end{proof}

\begin{lem}\label{strdiv} Let $R$ be an integral domain, $I$ a proper
strong ideal of $R$ and $x\in (R : I)\setminus I$.  Then $(I:_Rx)$ is also a
proper
strong ideal of $R$.  \end{lem}

\begin{proof} Clearly $J = (I:_Rx)$ is a proper ideal of $R$ containing
$I$. Hence $(R:J) \subseteq (R:I) =(I:I)$. Moreover,
$(R:J)Jx=(I:I)(I:_Rx)x\subseteq I$. Thus $(R:J)J\subseteq (I:_Rx)=J$,
which implies that $(R:J) =(J:J)$. \end{proof}

\begin{prop} \label{maxstrdiv} If $I$ is a maximal strongly divisorial
ideal of an integral domain $R$, then $I$ is a prime ideal of $R$.
\end{prop}

\begin{proof} If $I$ is a proper strongly divisorial ideal of $R$ and $x
\in R \setminus I$, then, by Lemma \ref{strdiv}, the ideal $(I:_Rx)$ is
proper and strongly divisorial.  We conclude by applying Lemma \ref{genprime}
to the set $\mathbf S$ of the proper strongly divisorial ideals of $R$.
\end{proof}

We observe that the proof of the previous proposition yields the
following more general:

\begin{prop}
For any $*$-operation, a maximal strong $*$-ideal is prime.
\end{prop}

We now give conditions sufficient for a strongly divisorial ideal to
have a strongly divisorial minimal prime.
It is known that  a prime ideal minimal over a
divisorial ideal is a $t$-ideal \cite{HK} and that a minimal prime of a
strong ideal is strong \cite[Proposition 2.1]{HKLM}.

\begin{lem}\label{divrad}
If a divisorial ideal $I$ of a domain $R$ contains a power of its radical,
then $\sqrt I$ is divisorial.
\end{lem}

\begin{proof}
By assumption, $(\sqrt I)^n\subseteq I$ for some $n\ge1$.  Thus $((\sqrt
I)_v)^n\subseteq ((\sqrt I)^n)_{v}\subseteq I$.  Hence $(\sqrt
I)_v\subseteq \sqrt  I$ and so $\sqrt I$ is divisorial.
\end{proof}

\begin{lem}\label{preirred}
If a divisorial ideal $I$ of a domain $R$ is
an intersection of primes and the prime $P$ is essential
for this representation (that is, it cannot be omitted), then $P$ is
divisorial.
\end{lem}
\begin{proof} 
 Let $I = \bigcap_{\lambda\in\Lambda}P_{\lambda}$
 be a divisorial ideal, where the $P_{\lambda}$ are primes,
 and let $\mu \in \Lambda$ be such that
$J:=\bigcap_{\{\lambda\ne\mu\}}P_\lambda \not \subseteq P_{\mu}$.
We have $(I:_{R}J)=\bigcap_{\{\lambda\in\Lambda\}}(P_\lambda:_{R}J)=
P_\mu$ since, for all $\lambda$, $J\subseteq P_\lambda$ if and only if
$\lambda\ne\mu$.
It follows that $P_\mu$ is
    divisorial.
 \end{proof}

\begin{lem}\label{irred}
    
    Let $\{P_{\lambda}\}_{\lambda\in\Lambda}$ be a family of prime ideals
and 
 let $I = \bigcap_{\lambda\in\Lambda}P_{\lambda}$
 be an irredundant intersection. Then the following conditions are
equivalent:
    \begin{enumerate}
    \item $I$ is strong (resp. divisorial);
    \item $P_{\lambda}$ is strong (resp. divisorial) for each
    $\lambda$.
    \end{enumerate}
    \end{lem}
    
    \begin{proof} 
        By \cite[Proposition 3.13]{HKLM}, $(R:I)$ is a ring if and
    only if $(R:P_{\lambda})$ is a ring for each
    $\lambda \in \Lambda$. But, if $J$ is a radical
    ideal, then $(R:J)$ is a ring if and
    only if $(R:J) = (J:J)$ \cite[Proposition 3.3]{A}. Hence $I$ is strong
if and 
    only if $P_\lambda$ is strong for each $\lambda$.
    
    If $P_{\lambda}$ is  divisorial for each
    $\lambda$, then clearly $I$ is divisorial. The converse follows
    from Lemma \ref{preirred}.
\end{proof}

Since the radical of a strong ideal is strong \cite[Proposition 2.1]{HKLM},
an immediate consequence of Lemmas \ref{divrad} and \ref{irred} is:

\begin{prop}\label{forte}
If a domain $R$ has a proper strongly divisorial ideal $I$
containing a power of its radical and this radical is an irredundant
intersection of primes, then all the minimal primes over $I$ are
strongly divisorial.
\end{prop}

We now prove that any seminormal domain satisfies the hypothesis that each
strongly
divisorial ideal is radical. Hence a seminormal domain which has a
strongly divisorial ideal that is an irredundant intersection of primes
has a strongly divisorial prime.

Recall that a domain $R$ is {\it seminormal} if
it contains all the elements  $x \in K$ such that $x^{2}, x^{3}
\in R$. 

\begin{lem}\label{rad}
 Let $R$ be a seminormal domain, $T$ an overring of $R$
 and $I = (R:T)$. Then $I$ is radical in both $T$ and $R$.
 \end{lem}
 \begin{proof} It is enough to show that $I$ is radical in $T$. Let
  $t$ be an element of $T$ such that $t^{2} \in I$, thus
 $t^{2}T \subseteq  R$, and
 $t^{3}T = t^{2}(tT) \subseteq  R$. Hence, for any $x
 \in T$, we have $(tx)^{2} = t^{2}x^{2} \in  R$ and
 $(tx)^{3}= t^{3}x^{3} \in  R$. Since $R$ is seminormal,
 it follows that $tx \in R$ and so $t \in I$.
\end{proof}
 
Lemma \ref{rad} leads to a characterization of
seminormal domains, which is similar to \cite[Theorem 1.1]{GH}.
 \begin{prop}\label{radical}
 A domain $R$ is seminormal if and only if  each
 strongly divisorial ideal $I$ of $R$ is radical in the overring $(I:I)$.
 \end{prop}
 \begin{proof}

If $I$ is a strongly divisorial ideal of $R$ and $T = (I:I)$, then $I =
(R:T)$. Conversely, if $T$ is an overring of $R$ so that $(R:T)\ne0$,
then the ideal $I = (R:T)$ is strongly divisorial and $T \subseteq
(I:I)$ \cite[Proposition 6]{B}.

Then, if $R$ is seminormal and $I$ is strongly divisorial,
$I$ is radical in $T = (I:I)$ by Lemma \ref{rad}.
 
 Conversely, assume that $R$ is not seminormal and let $x \in K
 \setminus R$ be such that $x^{2}, x^{3} \in R$. Then the ideal $I =
 (R: R[x])$ is not radical in $(I:I)$,
 because $x \in R[x] \subseteq (I:I)$, $x^{2} \in I$ but $x \notin I$.
 \end{proof} 

 We will give in the next section examples showing that both the
hypotheses of Proposition \ref{forte} are essential.
In fact the examples given in  \S2.1 and \S2.2 are non-completely
integrally closed, quasilocal one-dimensional domains, such that
 the maximal ideal $M$ is not divisorial. In these cases no proper
 strongly divisorial ideal contains a power of its radical $M$.

On the other hand the example given in \S2.3 is a non-completely integrally
closed  domain
 without strongly divisorial prime ideals which can be seminormal (or even
integrally closed).
In this case, no (radical) strongly divisorial ideal is an irredundant
intersection of prime ideals.

We now apply the previous results to obtain other classes of domains
that are completely integrally
closed if they have no strongly divisorial prime ideals.

\begin{thm}\label{noeth} Let $R\subsetneq T$ be integral domains with
$(R:T) \neq (0)$. If $T$ is Noetherian, then $R$ contains a strongly
divisorial prime ideal. \end{thm}

\begin{proof}
Let $I=(R:T)$, thus $I$ is a common proper ideal of $R$ and of $T$.
Moreover, $I$ is divisorial as an ideal of $R$.  Since $T$ is Noetherian,
the ideal $I$ of $T$ contains a power of its radical and it has
just finitely many minimal primes.  Since the radical of $I$ in $R$ is
contained in the radical of $I$ in $T$, the ideal $I$ of $R$ contains a
power of its radical and it is also an (irredundant)
intersection of finitely many primes of $R$.  By Proposition \ref{forte},
all  minimal primes over $I$ are strongly divisorial.
\end{proof}

We note that if the domains
$R$ and $T$ satisfy the hypotheses of the previous theorem, then
$R$ is  
a pullback of $T$. Moreover, $T\subseteq R^*$.




We now turn to the case of domains with all the maximal $t$-ideals
divisorial.

\begin{prop} \label{maxdiv} Let $R$ be a domain such that each maximal
$t$-ideal is divisorial.  If $R$ has no strongly divisorial prime ideals,
then $R$ is completely integrally closed.
\end{prop}

\begin{proof} Suppose, by contradiction, that $R$ is not completely
integrally closed, that is, that in $R$ there is a proper strongly
divisorial ideal.

Let $\{I_\lambda\}_{\lambda\in \Lambda}$ be a chain of proper strongly
divisorial ideals of $R$.  Then the ideal $I = \bigcup_{\lambda\in \Lambda}
I_\lambda$ is a $t$-ideal and, since each maximal $t$-ideal of $R$ is
divisorial, 
$I_v$ is a proper ideal.
Moreover, 
$$(R:I)=\bigcap_{\lambda\in \Lambda}(R:I_\lambda)=
\bigcap_{\lambda\in \Lambda}(I_\lambda:I_\lambda)\subseteq
\bigcap_{\lambda\in \Lambda}(I:I_\lambda)=(I:I).$$
Hence $I_v$ is a strongly divisorial ideal.

It follows that $R$ has a maximal strongly divisorial ideal, which is
prime by Proposition \ref{maxstrdiv}.
\end{proof}

\begin{thm} \label{Krull}

Let $R$ be any domain.  The following conditions are equivalent:

\begin{enumerate}
    \item $R$ has no proper strong $t$-ideals;

    \item $R$ has no strong prime $t$-ideals;

    \item $R$ has no strongly divisorial prime ideals and each maximal
$t$-ideal is $t$-invertible;

    \item $R$ has no strongly divisorial prime ideals and each  maximal
$t$-ideal is divisorial;

    \item$R$ is completely integrally closed and each maximal
    $t$-ideal is divisorial;
    
    \item $R$ is a Krull domain;
    
    \item Each
$t$-ideal of $R$ is $t$-invertible.  \end{enumerate} \end{thm}

\begin{proof}

(1) $\Rightarrow$ (2) and (3) $\Rightarrow$ (4) are clear.

(2) $\Rightarrow$ (3) follows from the facts that any  $v$-ideal is a
$t$-ideal and that a maximal $t$-ideal is either strong or $t$-invertible.

(4) $\Rightarrow$ (5) follows from the Proposition \ref{maxdiv}.

(5) $\Rightarrow$ (6) is \cite[Corollary 2.8]{G}.

(6) $\Rightarrow$ (7) is well known (see for example \cite [Corollary
2.7]{G}).

(7) $\Rightarrow$ (1) is clear.

\end{proof}

\smallskip 
The fact that a domain without strong prime $t$-ideals is completely
integrally
closed can be also proved recalling that a minimal prime of a proper
strongly divisorial ideal is a strong $t$-ideal \cite[Proposition
2.1]{HKLM}.

A domain in which each $t$-ideal is divisorial is called a
$TV$-domain.  Examples of $TV$-domains include Mori, in particular Noetherian,
domains and pseudovaluation domains \cite {HZ}.

The following corollary is immediate by Theorem \ref{Krull} and was
proved for Mori domains in \cite [Corollary 14]{B}; the equivalence (2)
$\Leftrightarrow$ (3) in it was also proved in \cite [Theorem 2.3]{HZ}.

\begin{cor} 
Let $R$ be any domain.  The following conditions are equivalent:

\begin{enumerate}
    \item $R$ is a $TV$-domain with no strongly divisorial prime ideals;

    \item $R$ is a completely integrally closed $TV$-domain;

    \item $R$ is a Krull domain.
\end{enumerate}

\end{cor}

An invertible maximal ideal is clearly a $t$-invertible maximal $t$-ideal.
Hence a domain in which each maximal ideal is invertible satisfies the
hypothesis that each maximal $t$-ideal is divisorial.  In addition, a
completely integrally closed domain in which each maximal ideal is
invertible is one-dimensional \cite[Theorem 2.1]{G}.  Hence we obtain the
following results.
    
\begin{cor} \label{Dedekind}

Let $R$ be any domain.  The following conditions are equivalent:

\begin{enumerate}
    \item $R$ has no strong ideals;

    \item $R$ has no strong prime ideals;

    \item $R$ has no strongly divisorial prime ideals and each maximal
    ideal is invertible;

    \item $R$ has no strongly divisorial prime ideals and each maximal
ideal is divisorial;

\item $R$ is completely integrally closed and each maximal ideal
is divisorial; 

\item $R$ is a Dedekind domain;

    \item Each ideal of $R$ is invertible.

 \end{enumerate}

\end{cor}

Recall that a $t$-ideal of a domain may be neither strong nor
$t$-invertible.  However, by Theorem \ref{Krull}, if no proper $t$-ideal is
strong, then all the $t$-ideals must be $t$-invertible.

\begin{cor}

Let $R$ be a domain in which each ideal is divisorial.  The following
conditions are equivalent:

\begin{enumerate}
    \item $R$ has no strongly divisorial prime ideals;

    \item $R$ is completely integrally closed;

    \item $R$ is a Dedekind domain.

 \end{enumerate}

\end{cor}




An integrally closed domain in which each ideal is divisorial is a Pr\"ufer
$\sharp\sharp$-domain \cite[Theorems 4.3.4 and 4.3.6]{FHP}.  Hence it has
the radical trace property \cite[Theorem 4.2.28]{FHP}.

Recall that $R$ is a $\sharp\sharp$-{\it domain} if each overring $S$ of
$R$ has the property that, for each pair of nonempty sets of maximal ideals
$\mone$ and $\mtwo$ of $S$, $\bigcap_{M \in \mone} S_{M} = \bigcap_{M \in
\mtwo} S_{M}$ only if $\mone=\mtwo$.

A domain satisfying the {\it radical trace property}, or an $RTP$-{\it
domain}, is a domain $R$ such that $I(R:I)$ is a radical ideal, for 
every
noninvertible ideal $I$.  That is, $R$ is an $RTP$-domain if and only if
each strong ideal of $R$ is radical \cite[Section 4.2]{FHP}.  The
property that $I(R:I)$ is a prime ideal for every noninvertible ideal $I$,
is called the {\it trace property} and domains with this property are
called $TP$-{\it domains}.  Therefore $R$ is a $TP$-domain if and only if
each proper strong ideal of $R$ is prime.  Examples of $TP$-domains 
include
valuation domains \cite[Section 4.2]{FHP}.

It follows directly from the definition that a $TP$-domain without strongly
divisorial primes is completely integrally closed.

We do not know whether in general an $RTP$-domain with no strongly
divisorial prime ideals is completely integrally closed.  However, the
answer to this problem is positive in case $R$ is an $RTP$-domain with the
ascending chain condition on radical ideals.  In fact, in this case, if the
set of maximal strongly divisorial  ideals of $R$ is not empty, it has
maximal elements and these are prime ideals by Proposition \ref{maxstrdiv}.
Recalling that each ideal of a domain with the
ascending chain condition on radical ideals has finitely many minimal
primes, the same result can also be obtained by applying Lemma
\ref{preirred}.
In addition,
it is known that, if $R$ is an $RTP$-domain, then each nonzero nonmaximal
prime of $R$ is strongly divisorial \cite [Theorem 4.2.16]{FHP}.  Hence an
$RTP$-domain without strongly divisorial primes is always
one-di\-men\-sional.
It follows that a Pr\"ufer $RTP$-domain without strongly divisorial primes
is completely integrally closed.




In general a one-dimensional domain without strongly divisorial prime
ideals may not be completely integrally closed, as we will show in the next
section. The following proposition is a partial positive result:

 
 \begin{prop}\label{seminormal}
  Let $R$ be a one-dimensional
  domain that is not completely integrally closed and
  such that the intersection of the maximal ideals is
  irredundant (e.g., $R$ has finitely many maximal ideals).
 If $R$ is either an $RTP$-domain or seminormal, then
  it has a strongly
  divisorial maximal ideal.
  \end{prop}
 \begin{proof} 
  If $R$ is not completely integrally closed,
  it has a strongly divisorial ideal
  which is radical (Lemma \ref{rad}). Since
  this ideal is an irredundant intersection of primes, we conclude by
  Lemma \ref{irred}.
  \end{proof}




\section{Counterexamples}

In this Section we give examples of domains which are not completely
integrally closed, but do not have strongly divisorial prime ideals.

We denote the integral closure of an integral domain $R$ in its quotient
field by $R'$.


\subsection{\hspace{-1.4cm}
\S \hspace{.8cm}A one-dimensional non-completely integrally closed
quasi\-local integral domain
$(R,M)$ such that $R^*=R'$ is also one-dimensional quasi\-local,
$(R:R^*) = (0)$ and such that $M$ is not divisorial.}

 \

Let $F$ be a field of finite characteristic $p$ and let $S$ be the additive
submonoid of $\mathbb Q$ generated by the set of the numbers
$\{n+1+\frac1{2^n}\}$
for all natural numbers $n\ge0$.
Clearly, $2,5 \in S$.  Thus each natural  number $n\ge4$ belongs to $S$ .

Consider the monoid domain $F[S]=F[X^S]$, where $X^S=\{X^s\,|\, s\in S\}$,
and set $M=X^{S\setminus\{0\}}F[S]$ and $R = F[S]_M$.

For $n\ge4$, since $n \in S$, then $X^n \in R$.  In particular this
shows directly
that $R$ is not seminormal, because $X \notin R$ (cf. Proposition
\ref{seminormal}).

In order to prove that $R$ has the required  properties, we need two Lemmas.

\begin{lem}\label{semigroup}

\
\begin{enumerate}
\item For $n\ge0$, if $m$ is a sufficiently large positive integer,
then $\frac m{2^n}\in S$.

\item Let $q$ be a rational number.  If $q+s \in S$ for
each nonzero $s \in S$, then $q \in S$.

\item
If $q$ is a rational number, then $q+\frac 1 {2^n}\notin S$ for
natural $n \gg 0$.
\end{enumerate} \end{lem}

\begin{proof}  

\begin{enumerate}
\item

For natural $k\gg0$, the elements $\frac {k2^n}{2^n}=k$ and
$\frac{k2^n+1}{2^n}=(k-n-1)+(n+1+\frac 1{2^n})$ belong to $S$ and the
numerators
of these two fractions are coprime positive integers.  Thus, if $m\gg0$,
then
$\frac m{2^n}\in S$.

\item

Let $q=\frac a b$, where $a$ and $b$ are integers, and $b>0$. Let
$k>\max(q,b)$ be an integer. Since $u:=q+k+1+\frac 1{2^k}\in S$, we have
$u=\sum_{i=1}^m(n_i+1+\frac 1{2^{n_i}})$ for some natural numbers
$n_1\le \dots \le n_m$, not necessarily distinct. Since $k>q$, we have
$n_i\ge k$ for at most one index $i$. Since $k>b$, we have $u=\frac a
b+k+1+\frac 1{2^k}=\frac \alpha{2^k\beta}$, where $\alpha$ and $\beta$
are odd integers. The denominator of $u=\sum_{i=1}^m(n_i+1+\frac
1{2^{n_i}})$ as a reduced fraction is at most $2^{n_m}$, thus $n_m\ge
k$, and $n_i<k$ for $i<m$. It follows that $n_m=k$ and that
$q=\sum_{i=1}^{m-1}(n_i+1+\frac 1{2^{n_i}})\in S$.

\item
 Let $q=\frac a b$,
where $a$ and $b$ are integers, and $b>0$.
Let $n>\max(q,b)$ be an integer.
  Assume that $q+\frac 1{2^n}\in S$, thus $q+\frac 1{2^n}=
\sum_{i=1}^m(n_i+1+\frac 1{2^{n_i}})$
for some natural numbers $n_1\le\dots\le n_m$.
Since $2^n$ divides the denominator of $q+\frac 1{2^n}$ written as a
reduced fraction, we see that $n_m\ge n$.
 Thus $n_m+1+\frac {1}{2^{n_m}}> n+ \frac {1}{2^n}>
q+\frac {1}{2^n}$, a contradiction.

\end{enumerate}
\end{proof}

\begin{lem}\label{ass}
Let $f$ be a nonzero element of $R$.  Then, for $n$ sufficiently large,
$f^{p^n}$ is associated in $R$ with an element of the form $X^s$ for some
$s\in S$.

\end{lem}

\begin{proof}
We may assume that $f$ belongs to the ideal $M$ of $F[S]$, so that
$$f=a_1X^{s_1}+a_2X^{s_2}+\dots +a_mX^{s_m},$$ where $s_1<\dots<s_m$ are
elements of $S$, and the coefficients $a_1,\dots,a_m$ are in $F$. We may
further assume that $a_1=1$.  We have
$$f^{p^n}=X^{p^ns_1}+a_2^{p^n}X^{p^ns_2}+\dots +a_m^{p^n}X^{p^ns_m}$$ $$=
X^{p^ns_1}(1+a_2^{p^n}X^{p^n(s_2-s_1)} +\dots
+a_m^{p^n}X^{p^n(s_m-s_1)}).$$ If $n$ is sufficiently large, by Lemma
\ref{semigroup} (1), we see that $p^n(s_i-s_1) \in S$, for $i=2, \cdots, m$
and so $$1+a_2^{p^n}X^{p^n(s_2-s_1)} +\dots +a_m^{p^n}X^{p^n(s_m-s_1)}\in
F[S].$$ Hence $f^{p^n}$ is associated in $R$ with $X^{p^ns_1}$.
\end{proof}

\begin{prop}

\

\begin{enumerate}
\item $R$ is one-dimensional.

\item $R'=R^*$ is one-dimensional quasilocal.

\item $(R:R^*) = (0)$.

\item The ideal $M$ is not divisorial.

\item $R$ is not integrally closed and has no strongly
divisorial prime ideals.  \end{enumerate} \end {prop}

\begin{proof} 

\begin{enumerate}
\item

If $f,g$ are nonzero elements of $M$, then by Lemma \ref{ass} there exist
$m,n \in \mathbb N$ such that $f^m$ and $g^n$ are associated.  So $\sqrt
{(f)}$ = $\sqrt {(g)}$ and $R$ is one-dimensional (cf. \cite[Theorems 17.1
and 21.4]{Gcsg}).

\item

By Lemma \ref{semigroup} (1) and by \cite[Corollary 12.7]{Gcsg}, we
easily obtain that $R'=R^*=F[X^{\frac m n}\ |\ \text{$m$ and $n$ are
positive integers}]_N$, where $N$ is the maximal ideal generated by the
elements $X^{\frac m n}$ for positive integers $m$ and $n$. Hence
$R^*=R'$ is one-dimensional quasilocal.

\item Suppose that $0 \neq f = a_1X^{s_1}+a_2X^{s_2}+\dots +a_mX^{s_m}
\in (R:R^*)$, where $a_1\ne0$. Then, for $n\gg0$, $X^{\frac 1{2^n}} \in
R^*$, but $fX^{\frac 1{2^n}} =
a_1X^{s_1+\frac 1{2^n}}+ \cdots \notin R$, since $s_1+\frac 1{2^n}\notin
S$ by Lemma \ref{semigroup} (3); a contradiction.

\item

 To show that $M$ is not divisorial, it is enough to prove that $(R:M)
 \subseteq R$.

 Let $t=\frac f g$, where $f,g\in F[S]$, such that $tM \subseteq R$.  By
 Lemma \ref{ass}, for $n\gg0$, the element $g^{p^n}$ is associated in $R$
 with an element of $X^S$.  Since $t=\frac{fg^{p^n-1}}{g^{p^n}}$, we may
 assume that $t\in F[G]$, where $G$ is the group of quotients of $S$.  Since
 $M=X^{S\setminus\{0\}}R$, it follows that $t^\prime M\subseteq M$ for all
components
 $t^\prime\in X^G$ of $t$.  Thus we may assume that $t=X^q$, where $q\in G$,
 thus $q$ is a rational number.  Since $tX^s \in R$, for each nonzero $s \in S$,
 by Lemma \ref {semigroup} (2) we obtain $t \in R$.

\item

Since, as we have already observed, $R$ is not seminormal, it is 
certainly not integrally closed.  Since $R$ is one dimensional 
and quasilocal and
the unique maximal ideal $M$ is not divisorial, there are no strongly
divisorial prime ideals.

\end{enumerate}

\end{proof}


\subsection{\hspace{-1.4cm}
\S \hspace{.8cm}A one-dimensional non-completely integrally closed
quasi\-local domain $(R,M)$ such
that $R^*$ is also one-dimensional quasi\-local,
$(R:R^*)\ne (0)$, and such that $M$ is not divisorial.}
\

Let $\mathbf X=\{X_n,|\,n\ge1\}$ be a set of indeterminates over a field
$F$.  We define $w(X_n)=\frac 1 n$ (the {\em weight} of $X_n$) for all
$n\ge1$.  We extend the function $w$ to the set of monomials in $F[\mathbf
X]$ setting $w(cX_{i_1}\dots X_{i_m})=\sum_{j=1}^m w(X_{i_j})$, where $c$ is
a nonzero element of $F$ and $X_{i_1}\dots X_{i_m}$ are indeterminates in
the set
$\mathbf X$ ($i_1,\dots, i_m$ not necessarily distinct).  For a nonzero
polynomial $ f\in F[\mathbf X]$, we define
$w(f)$ to be the least weight of a monomial occurring in $f$.  We let
$w(0)=\infty$. Clearly, $w(fg)=w(f)+w(g)$ for any polynomials $f$ and
$g$.  We uniquely extend the function $w$ from $F[\mathbf X]\setminus
\{0\}$ to a function defined on $F(\mathbf X)$ and satisfying the above
property for rational functions $f$ and $g$.  Using the weight $w$ we
naturally define a $\mathbb Q^+$-grading on $F[\mathbf X]$: $$ F[\mathbf
X]=\bigoplus_{q\in\mathbb Q^+}F[\mathbf X]_q, $$ where $F[\mathbf X]_q$ is
the set of {\em
$w$-homogeneous} forms of weight $q$, together with $0$, that is, the
$F$-linear combinations of monomials of weight $q$ for $q\in\mathbb Q^+$.
If $f=\sum_{i=0}^rT_{q_i}\in F[\mathbf X]$, where $0\le q_0<q_1<\dots<q_r$
are rational numbers, and $T_{q_0},\dots,T_{q_r}$ are $w$-homogeneous forms
of weight $q_0,\dots,q_r$ respectively, we set $\win(f)=T_{q_0}$ ({\em the
$w$-initial form} of $f$).  We let $\win(0)=0$.  Clearly,
$\win(fg)=\win(f) \win(g)$ for any polynomials $f$ and $g$.  Thus we may
define $\win(\frac f g)=\frac{\win(f)}{\win(g)}$ for any nonzero rational
function $\frac f g$, where $f$ and $g$ are nonzero polynomials.

If $q\in \mathbb Q$ , we denote by $F(\mathbf X)_{\ge q}$ the set $\{f\in
F(\mathbf X)|\ w(f)\ge q\}$.

Let $A=F[\mathbf X]+F(\mathbf X)_{\ge 1}$, and let $P$ be the maximal ideal
$(\mathbf X)F[\mathbf X]+F(\mathbf X)_{\ge 1}$ of $A$.  Set $R=A_P$, and
$M=PA_P$, thus $(R,M)$ is a quasilocal domain with quotient field
$F(\mathbf X)$.  Moreover:

\begin{prop}

\ 
\begin{enumerate}

\item $R^*=F(\mathbf X)_{\ge 0}$ and $(R:R^*) \neq (0)$.

\item  $R$ and $R^*$ are quasilocal and
one-dimensional, and $R^*$ dominates $R$.

\item The ideal $M$ is not divisorial.

\item $R$ is not completely integrally closed and has no strongly
divisorial prime ideals.  \end{enumerate} \end {prop}

\begin{proof}

\begin{enumerate}
\item

Clearly $R\subseteq F(\mathbf X)_{\ge 0}$.  The domain $F(\mathbf X)_{\ge
0}$ is completely integrally closed.  Indeed, let $0\ne f\in F(\mathbf X)$.
If $gf^n\in F(\mathbf X)_{\ge 0}$ for some nonzero $g\in F(\mathbf X)_{\ge
0}$ and all $n\ge1$, then $w(gf^n)=w(g)+nw(f)\ge0$ for all $n\ge1$, which
implies that $w(f)\ge0$.  It follows that $f\in F(\mathbf X)_{\ge 0}$.

Since $X_1 F(\mathbf X)_{\ge 0}\subseteq R$, then $F(\mathbf X)_{\ge
0}\subseteq R^*$.  Hence $R^*=F(\mathbf X)_{\ge 0}$, and $(R:R^*)\ne0$.

\item
Clearly, $N:=F(\mathbf X)_{>0}$ is the unique maximal ideal of $F(\mathbf
X)_{\ge 0}$, and $M:=N\cap R$ is the unique maximal ideal of $R$.  If $c_1$
and $c_2$ are nonzero elements in $N=F(\mathbf X)_{>0}$, then, for
sufficiently large $n$, we have $w(\frac{c_1^n}{c_2})\ge1$.  Hence
$\frac{c_1^n}{c_2}\in R$.  This implies that both $R$ and $F(\mathbf
X)_{\ge 0}$ are one-dimensional since any element in $R\setminus F(\mathbf
X)_{>0}$ is a unit in $R$.

\item Assume that the maximal ideal $M$ of $R$ is divisorial, i.e. $(R:M)
\neq R$.

Let $h$ be an element in $F(\mathbf X)\setminus R$ satisfying
$hM\subseteq R$. First we show that $\win(h)\in F[\mathbf X]$. Otherwise we have
$w(h)<1$. Let $m$ be a positive integer such that $w(h)+\frac 1m<1$.
Choose $n\ge m$ so that $h\in F(X_1,\dots,X_{n-1})$. Since $X_n\in M$,
we have $hX_n\in M$, and so there are an element $u\in 1+(\mathbf
X)F[\mathbf X]$ and a polynomial $p\in F[\mathbf X]$ such that
$uhX_n-p\in F(\mathbf X)_{\ge1}$. Since $\win(hX_n)<1$, we obtain that
$\win(p)=\win(uhX_n)=X_n\win(h)$. Since $h\in F(X_1,\dots,X_{n-1})$, we
infer that $\win(h)\in F[\mathbf X]$, a contradiction. Hence $\win(h)\in F[\mathbf X]$ as claimed.
For some $m$ we have $h\in F(X_1,\dots,X_m)$.
Define inductively $h_0=h$ and $h_{n+1}=h_n-\win(h_n)$. Clearly,
$h_n\in F(X_1,\dots,X_m)\setminus R$, $h_nM\subseteq R$, and $w(h_{n+1})>w(h_n)$ for all
$n\ge0$.  It follows that
$\win(h_n)\ge1$, so $h_n\in R$ for $n\gg0$, a contradiction.

\item By part (1), $R$ is not completely integrally
closed; for example, $\frac {X_1}{X_2^2}\in R^*\setminus R$. Further, as
in the previous example, since $R$ is one dimensional quasilocal and the
unique maximal ideal $M$ is not divisorial, there are no strongly
divisorial prime ideals.

\end{enumerate}

\end{proof}


\subsection{\hspace{-1.4cm}
\S \hspace{.8cm}An infinite dimensional integrally closed
non-completely integrally closed domain without 
strongly divisorial prime ideals
such that $(R:R^*) \neq (0)$.}

\

To construct this example, we use some properties of the domain $E = D[Z,
\frac dZ]$, where $Z$ is an indeterminate over $D$ and $d \in D$ is a
nonzero 
element. This is isomorphic to the extended Rees Algebra $D[dX,\frac 1X]$ of
the
nonzero principal ideal $dD$, via the correspondence $Z \to \frac 1X$.

We have $\Qf(E) = (\Qf(D))(Z)$ and $D = E \cap \Qf(D)$.

\smallskip

If $(D:D^*) \neq (0)$, then $E^* = D^*[Z,
\frac dZ]$ and $(D:D^*) \subseteq (E:E^*)$.
In fact, since $D^*$ is completely integrally closed,
then $D^*[Z, \frac dZ]$ is completely integrally closed \cite
[Theorem 8]{AA} and $D^*[Z, \frac dZ] \subseteq
E^*$.

\medskip

The domain $R$ is constructed in the following way.

\medskip

Let $A$ be an integral domain
such that $A \neq A^*$ and  $(A:A^*) \neq (0)$, and let $a$ be a nonzero
element of $(A:A^*)$.

Let $T_0$ and $T_{\varepsilon_1\dots\varepsilon_n 0}$, for $n\geq 1$ and
$\varepsilon_1,\dots,\varepsilon_n \in \{0,1\}$, be independent
indeterminates over $A$.

\smallskip

Define inductively $T_1 = \frac {a}{T_0}$, and
$T_{\varepsilon_1\dots\varepsilon_n1}= \frac
{T_{\varepsilon_1\dots\varepsilon_n}}{T_{\varepsilon_1\dots\varepsilon_n
0}}$ for $n \geq 1$ and set $$R =
A[T_{\varepsilon_1\dots\varepsilon_n}\,|\,n\ge1 \text{ and }
\varepsilon_1,\dots,\varepsilon_n \in \{0,1\}].$$

\medskip
 We have $R=\bigcup_{n=0}^\infty R_n$, where $R_0 = A$ and
$$R_n= A[T_{\varepsilon_1\dots\varepsilon_k}\,|\,1\leq k \leq n \text{ and
} \varepsilon_1,\dots,\varepsilon_k \in \{0,1\}]$$ for each $n\geq 1$.

Moreover the domain $R_n$ is obtained from $A$ by constructing iteratively
extended Rees algebras of the type $D[Z, \frac dZ]$; in fact,
$$R_{1}=
A[T_0,T_{1}] = A[T_0,\frac {a}{T_0}]$$
and, for $n \geq 1$,
$$R_{n+1}=R_n[Z_1,\dots,Z_{2^n},\frac
{c_1}{Z_1},\dots,\frac{c_{2^n}}{Z_{2^n}}], $$
where 
$$\{Z_1,\dots,Z_{2^n}\}
= \{T_{\varepsilon_1\dots\varepsilon_{n}0}|\,
\varepsilon_1,\dots,\varepsilon_{n} \in \{0,1\}\}$$
is a set of
algebraically independent elements over $R_n$, and
$$\{c_1,\dots,
c_{2^n}\} = \{T_{\varepsilon_1\dots\varepsilon_{n}}|\,
\varepsilon_1,\dots,\varepsilon_{n} \in \{0,1\}\} \subseteq R_n ,$$
because by definition $T_{\varepsilon_1\dots\varepsilon_{n}1}= \frac
{T_{\varepsilon_1\dots\varepsilon_{n}}}
{T_{\varepsilon_1\dots\varepsilon_{n} 0}}$.

\medskip

Setting $K_{0}= \Qf(A)$ and $K_{n} = \Qf(R_{n})$ for $n \geq
1$ , we have $K_{n+1} = K_{n}(Z_1,\dots,Z_{2^n})$ and $\Qf(R) =
\bigcup_{n=0}^\infty K_n.$
Since $K_{n} \cap R_{n+1}=R_{n}$, it follows that $K_{n} \cap
R = R_{n}$, for all $n \geq 0$.

\medskip

With this notation, we have:

\begin {prop}

\

\begin{enumerate}

\item $R^*=\bigcup_{n=0}^\infty R_n^*$.

\item $R$ is not completely integrally closed and $a \in (R :R^*)$.
In particular $(R :R^*)\neq (0)$.

\item $R$ has no strongly  divisorial prime ideals.

 \end{enumerate}

\end{prop}

\begin{proof}

\begin{enumerate}
\item
Let $x \in \Qf(R)$ be almost integral over $R$ and let $d \in R$ such that
$dx^k \in R$, for all $k \geq 0$.  Since $x \in \Qf(R_{s})$ and $d \in
R_{t}$ for some $s, t \geq 0$, then, for $m$ sufficiently large, $x \in
\Qf(R_m)$ , $d \in R_m$ and $dx^k \in R \cap \Qf(R_m) = R_{m}$, for all $k
\geq 0$.  Hence $x \in R_m^*$ and $R^*\subseteq \bigcup_{n=0}^\infty
R_n^*$.

The opposite inclusion is clear.

\item If $R = R^*$, then $A^{*} \subseteq R^* \cap \Qf(A) = R \cap \Qf(A)
= A$ .  Since $A$ is not completely integrally closed, neither is $R$.

By hypothesis, $a \in (A :A^*)$ and we have $$(A :A^*) \subseteq (R_1
:R_1^*) \subseteq (R_2 :R_2^*) \subseteq \dots $$ since
$R_n$ is obtained from $A$ by constructing iteratively
extended Rees algebras of the type $D[Z, \frac dZ]$.

To conclude it is enough to observe that $$\bigcup_{n=0}^\infty (R_n
:R_n^*) \subseteq \bigcup_{n=0}^\infty (R :R_n^*) \subseteq (R :R^*).$$
This follows from  (a): $R^*=\bigcup_{n=0}^\infty
R_n^*$.

Actually, $(R :R^*) = \bigcup_{n=0}^\infty (R_n :R_n^*)$.  Indeed, let $y
\in (R :R^*)$. If $y \in \Qf(R_{m})$ , then $yR_{m}^* \subseteq R \cap
\Qf(R_{m}) = R_{m}$ .  Hence $y \in (R_{m}:R_{m}^*)$ and it follows that
$(R :R^*) \subseteq \bigcup_{n=0}^\infty (R_n :R_n^*)$.

\item

Let $P$ be a strongly divisorial prime ideal of $R$.  Since $P$ contains
$(R :R^*)$, by part (2), $a\in P$.

Since $a=T_0T_1$, either $T_0$ or $T_1$ belongs to $P$.  We may assume that
$T_0\in P$ since there is an automorphism of $R$ over $A$ interchanging
$T_0$ and $T_1$.  Then, since $T_0=T_{00}T_{01} \in P$, either $T_{00}$ or
$T_{01}$ belongs to $P$.  Since there is an automorphism of $R$ over
$A[T_0]$ interchanging $T_{00}$ and $T_{01}$, we may assume that $T_{00}
\in P$.  Iterating this process, we may assume that all the indeterminates
$T_0,T_{00},T_{000}\dots$ are in $P$.  Thus $(R:P) \subseteq
(R:\{T_0,T_{00},\dots\})$.

We now prove that $(R:\{T_0,T_{00},\dots\}) = R$, thus obtaining a
contradiction.

Let $f\in (R:\{T_0,T_{00},\dots\})$ and let $n\ge 0$ such that $f\in
\Qf(R_n)$.
As above, we can write
$$R_{n+1}=R_{n}[Z_1,\dots,Z_{2^n},\frac
{c_1}{Z_1},\dots,\frac{c_{2^n}}{Z_{2^n}}],$$ where we can assume that
$Z_1=T_{0\dots 0}$.

Hence $fZ_1\in R \cap \Qf(R_{n+1}) = R_{n+1}$.  Since the indeterminates
$Z_1,\dots,Z_{2^n}$ do not occur in $f$, setting all the $Z_i$'s equal to
$1$, we obtain that $f\in R_n\subseteq R$.  \end{enumerate}

\end{proof}

To finish, we show in Proposition \ref{semarch} below that this domain
can be chosen to satisfy or not either one of the conditions of being
seminormal, integrally closed or Archimedean. In the seminormal case the conductor $(R:R^*)$
is a radical strongly divisorial ideal of $R$ (Lemma \ref{rad}), that is not
an
 irredundant intersection of primes (Lemma
\ref{irred}).

\begin{lem}\label{reessemarch}

Let $D$ be an integral domain, $d\in D\setminus\{0\}$ and $E = D[Z,\frac
dZ]$.  Then 
\begin{enumerate}
\item  $E$ is seminormal if and only if $D$ is
seminormal.  
\item $E$ is Archimedean if and only if $D$ is
Archimedean.  

\end{enumerate}
\end{lem}

\begin{proof} 
 \begin{enumerate}
\item Since $D = E \cap Qf(D)$, if $E$ is seminormal, then $D$ is
seminormal.

Conversely, assume that $D$ is seminormal and let  $f \in Qf(E)$ be such
that $f^{2},f^{3} \in E$. Since $E \subseteq D[Z, \frac 1Z]$ and
$D[Z, \frac 1Z]$ is seminormal, then $f \in D[Z, \frac 1Z]$. Write
$f = \sum_{i = m}^n c_{i}Z^{i}$. Since $D[Z] \subseteq E$,
we can assume that $m < 0$. We have that $f^{2} = c_{m}^{2}Z^{2m} + g
\in E$. Hence $c_{m}^{2} \in Dd^{|2m|}$ and similarly
$c_{m}^{3} \in Dd^{|3m|}$. Since $D$ is seminormal, we obtain that
$\frac {c_{m}}{d^{|m|}} \in D$ and so $c_{m}Z^{m} \in E$.

Repeating this process, we get that all the terms of $f$ are in $E$
and so $f \in E$.

\item Since $D = E \cap Qf(D)$, if $E$ is Archimedean, then $D$ is
Ar\-chi\-me\-dean.

Conversely, assume that $D$ is Archimedean. Let
$f(Z),g(Z)$ be two nonzero Laurent polynomials in $E$ such that
$g\in\bigcap_{n=1}^\infty Ef^n$. Since the domain $\Qf(D)[Z,\frac 1
Z]$ is Notherian,  so Archimedean, we see that $f=cZ^m$ for some element
$c\in D$ and some integer $m$. Since the coefficients of $g$ belong to
$\bigcap_{n=1}^\infty Dc^n$, we see that $c$ is invertible in $D$. If
$m\le0$, then clearly $f$ is a unit in $E$. If $m>0$, then the
coefficients of $g$ belong to $\bigcap_{n=1}^\infty Dd^n$. Hence $d$ is a
unit in $D$, $Z$ is a unit in $E$, and $f=cZ^m$ is a unit in $E$. We
conclude that $E$ is Archimedean.

\end{enumerate}
\end{proof}

\begin {prop}\label{semarch}
\
\begin{enumerate}
\item $R$ is seminormal (resp.  integrally closed) if and only if $A$ is
seminormal (resp.  integrally closed).

\item $R$ is Archimedean if and only if $A$ is Archimedean.

\end{enumerate}

\end{prop}

\begin{proof}

\begin{enumerate}
\item

Since $R =\bigcup_{n=0}^\infty R_n$, it is enough to show that $R_n$ is
seminormal (resp.  integrally closed) if and only if $R_{n+1}$ is
seminormal (resp.  integrally
closed), for each $n\ge 0$ .

This follows from the fact that the domain $R_{n+1}$ is
obtained from $R_n$ by constructing iteratively extended Rees algebras of
principal ideals and from the fact that $D[Z, \frac dZ]$ is seminormal
(resp.  integrally closed) if and only if $D$ is seminormal (resp.
integrally closed) by Lemma \ref{reessemarch} (resp. \cite[Theorem
8]{AA}).

\item
Since $A = R \cap Qf(A)$ , if $R$ is Archimedean, then
$A$ is Archi\-medean.

Conversely, assume that $A$ is Archimedean.  By Lemma \ref{reessemarch}, all
the domains $R_m$ are Archi\-me\-dean.

Now let $r$ be a nonzero element of $R$, and let $c\in \bigcap_{n=1}^\infty
Rr^n$.  For some $m\ge0$ we have $r,c\in R_m$.  Since $\Qf(R_m)\cap R=R_m$,
we see that $c\in \bigcap_{n=1}^\infty R_mr^n$.  Since $R_m$ is
Archimedean, we conclude that $c=0$ and that $R$ is Archimedean.
\end{enumerate}

\end{proof}



\end{document}